\documentclass[11pt]{article}
\usepackage{amsfonts}
\usepackage{latexsym,amscd,amsmath,amssymb}
\title{Sharp inequalities of homogeneous expansions for quasi-convex mappings of type B and almost starlike mappings of order $\alpha$}
\author{Ming-Sheng Liu$^{1}$\thanks{Corresponding author. E-mail address: liumsh@scnu.edu.com. }, Fen Wu$^{1}$  and Yan Yang$^{2}$\\
 {\it 1. School of Mathematical Sciences, South China Normal } \\
 {\it  University, Guangzhou 510631, China}\\
 {\it 2. School of Mathematics and Computational Science, } \\
 {\it Sun Yat-Sen University, Guangzhou 510275, China}
 }

\date{}
\begin{document}
\maketitle
\def\leq{\leqslant}
\def\geq{\geqslant}
\def\no{\noindent}
\def\dl{\displaystyle\lim}
\def\ds{\displaystyle\sum}
\def\di{\displaystyle\int}
\newcommand{\f}{\displaystyle\frac}
\newcommand{\il}{\int\limits}
\newcommand{\s}{\sum\limits}
\newcommand{\ba}{\begin{array}}
\newcommand{\ea}{\end{array}}
\newcommand{\ld}{\left\{}
\newcommand{\rd}{\right\}}
\newcommand{\lz}{\left[}
\newcommand{\rz}{\right]}
\newcommand{\lx}{\left(}
\newcommand{\rx}{\right)}
\newcommand{\og}{\omega}
\newcommand{\ra}{\rightarrow}
\newcommand{\he}[2]{\sum\limits_{#1}^{#2}}
\newcommand{\fd}[4]{\left \{ \begin{array}{ll}
#1&#2\\#3&#4\end{array}\right.}
\renewcommand{\thefootnote}{\fnsymbol{footnote}}

\footnote[0]{This research is supported by Guangdong Natural Science Foundation (Grant No.2014A030307016, 2014A030313422).}

\vspace{0.5cm}
\begin{center}
\begin{minipage}{15cm}
{\bf Abstract}\ \quad In this paper, we first obtain several sharp inequalities of homogeneous expansion for both the subclass of all normalized biholomorphic quasi-convex mappings of type B and order $\alpha$ and the subclass of all normalized biholomorphic almost starlike mappings of order $\alpha$ defined on the unit ball $B$ of a complex Banach space $X$. Then, with these sharp inequalities, we derive the sharp estimates of the third and fourth homogeneous expansions for the above mappings defined on the unit polydisk $D^n$ in $\mathbb{C}^n$.



\end{minipage}
\end{center}

\section{Introduction}
\qquad  It is well-known that there exists the following {\it de  Branges theorem} (or {\it Bieberbach conjecture}) in geometric function theory of one complex variable.

{\bf Theorem~A\cite{GK2003}}\quad If $f(z)=z+\sum_{n=2}^{\infty}a_{n}z^n$ is a biholomorphic function on the unit disk
$D=\{z\in \mathbb{C}: |z|<1\}$, then
\begin{eqnarray*}
|a_{n}|\leqslant n,\quad n=2,3,\cdots.
\end{eqnarray*}
The above estimations are sharp with the Koebe function $K(z)=\frac{z}{(1-z)^2}=z+\sum_{n=2}^{\infty}n\, z^n$ as their extremal function.

Unfortunately, H. Cartan\cite{CH1933} has pointed out that the above theorem does not hold in the case of several complex variables. Therefore, it is necessary to require some additional properties of mappings of a class in order to obtain the analogous estimations, for instance, the starlikeness, the
convexity and so on. In 1998, Sheng Gong posed the following conjecture.

{\bf Conjecture~B}\quad If $f:D^n\rightarrow\mathbb{C}^n$ is a normalized biholomorphic starlike mapping, where
$D^n=\{z=(z_1,\cdots, z_n)\in\mathbb{C}^n:\, |z_k|<1,\, k=1, 2, \cdots, n\}$ is the unit poly-disk in $\mathbb{C}^n$, then
\begin{eqnarray*}
\frac{\|D^{m}f(0)(z^m)\|}{m!}\le m\|z\|^m, z\in D^n, m=2,3,\cdots.
\end{eqnarray*}

We shall refer to Conjecture B as {\it Bieberbach-Gong Sheng Conjecture}. At the present time, only the case of $m = 2$ (see \cite{G1999}) has been shown. However, with respect to the estimation of homogeneous expansion for normalized biholomorphic starlike mappings on the Euclidean unit ball $B^n$ in $\mathbb{C}^n$, Roper and Suffridge in \cite{RS1999} have provided a counter example to verify that the above conjecture does not hold for $m = 2$. Then, in this paper, we extend some results in one complex variable to the case in several complex variables on the unit polydisk $D^n$ in $\mathbb{C}^n$. From these results, we also obtain the sharp third estimates of homogeneous expansion for biholomorphic quasi-convex mappings of type B and order $\alpha$, and the sharp third and fourth estimates of homogeneous expansion for biholomorphic almost starlike mappings of order $\alpha$ on the unit polydisk $D^n$ in $\mathbb{C}^n$.

We first recall the following results in the case of one complex variable.

{\bf Theorem~C}\cite{GK2003}\quad If $f(z)=z+\sum_{n=2}^{\infty}a_{n}z^n$ is a univalent function on the unit disk $D$ in the complex plane $\mathbb{C}$, then
\begin{eqnarray*}
|a_3-a_{2}^2|\le 1.
\end{eqnarray*}

{\bf Theorem~D}\cite{GK2003}\quad If $f(z)=z+\sum_{n=2}^{\infty}a_{n}z^n$ is a univalent convex function on the unit disk $D$, then
\begin{eqnarray*}
|a_3-a_{2}^2|\le\frac{1-|a_2|^2}{3}.
\end{eqnarray*}

{\bf Corollary~E}\cite{LL2013}\quad If $f(z)=z+\sum_{n=2}^{\infty}a_{n}z^n$ is a univalent convex function on the unit disk $D$ in $\mathbb{C}$, then
\begin{eqnarray*}
|a_3-\frac{2}{3}a_{2}^2|\le \frac{1}{3}.
\end{eqnarray*}

{\bf Corollary~F}\cite{LL2013}\quad If $f(z)=z+\sum_{n=2}^{\infty}a_{n}z^n$ is a univalent starlike function on the unit disk $D$ in $\mathbb{C}$, then
\begin{eqnarray*}
|2a_3-a_{2}^2|\le2
\end{eqnarray*}

{\bf Theorem~G}\cite{GK2003}\quad If $f(z)=z+\sum_{n=2}^{\infty}a_{n}z^n$ is a univalent convex function of order $\alpha$ on the unit disk $D$ in $\mathbb{C}$, $0\le\alpha<1$, then
\begin{eqnarray*}
|a_n|\le\frac{1}{n!}\Pi_{k=2}^{n}(k-2\alpha), n=2,3,\cdots.
\end{eqnarray*}
The above estimation is sharp.

{\bf Theorem~H}\cite{GK2003}\quad If $f(z)=z+\sum_{n=2}^{\infty}a_{n}z^n$ is a univalent starlike function of order $\alpha$ on the unit disk $D$ in $\mathbb{C}$, $0\le\alpha<1$, then
\begin{eqnarray*}
|a_n|\le\frac{1}{(n-1)!}\Pi_{k=2}^{n}(k-2\alpha), n=2,3,\cdots.
\end{eqnarray*}
The above estimation is sharp.

In the following, let $X$ be a complex Banach space with norm $\|\cdot\|$, $B=\{x\in X: \|x\|<1\}$ be the open unit ball in $X$; let $\partial_0B$ be the boundary of $B$, $\bar{B}$ be the closure of $B$, and $\partial_0D^n$ be the characteristic boundary (i.e. the boundary
on which the maximum modulus of the holomorphic function can be attained) of the unit poly-disk $D^n$. $\mathbb{N}_+$ represents the set of positive integers. Let
$H(B)$ be the set of all holomorphic mappings from $B$ into $X$, It is well known that if $f\in H(B)$, then
\begin{eqnarray*}
f(y)=\sum_{n=0}^{\infty}\frac{1}{n!}D^nf(x)((y-x)^n)=\sum_{n=0}^{\infty}\frac{1}{n!}D^nf(x)(\underbrace{y-x,\cdots,y-x}_{n})
\end{eqnarray*}
for all y in the neighborhood of $x\in B$, where $D^nf(x)$ is the $nth-Fr\acute{e}chet$ derivative of $f$ at $x$. Moreover, $D^nf(x)$ is a bounded symmetric $n-$linear mapping from $\prod_{j=1}^{n}X$ into $X$.

A holomorphic mapping $f:B\to X$ is said to be biholomorphic if the inverse $f^{-1}$
exists and is holomorphic on the open set $f(B)$. A mapping $f\in H(B)$ is said to be
locally biholomorphic if the $Fr\acute{e}chet$ derivative $Df(x)$ has a bounded inverse for each
$x\in B$. If $f : B\to X$ is a holomorphic mapping, we say that $f$ is normalized if
$f(0) = 0$ and $Df(0)=I$, where I represents the identity operator from $X$ into $X$.

If $X^*$ is the dual space of $X$, for each $x\in X\backslash\{0\}$, we define
\begin{eqnarray*}
T(x)=\{T_x\in X^*:\|T_x\|=1, T_x(x)=\|x\|\}.
\end{eqnarray*}
According to Hahn-Banach theorem, $T(x)$ is nonempty. For any $\alpha(\ne0)\in \mathbb{C}$, since $\frac{|\alpha|}{\alpha}T_x\in T(\alpha x)$ corresponding to each $T_x\in T(x)$, we always denote $\frac{|\alpha|}{\alpha}T_x$ by $T_{\alpha x}$.

Let $\Omega\subset\mathbb{C}^n$ be a bounded circular domain. The first order Fr$\acute{e}$chet derivative and the $m(m\ge2)$ order Fr$\acute{e}$chet derivatives of the mapping $f\in H(\Omega)$ were denoted by $Df(z)$ and $D^mf(z)(b^{m-1},\cdot)$, respectively. The corresponding matrixes become
\begin{eqnarray*}
&&Df(z)=\Big{(}\frac{\partial f_p(z)}{\partial z_k}\Big{)}_{1\le p,k\le n},\\
&&D^mf(z)(a^{m-1},\cdot)=\Big{(}\sum_{l_1.l_2,\cdots,l_{m-1}}^n\frac{\partial^mf_p(z)}{\partial z_k\partial z_{l_1}\cdots\partial z_{l_{m-1}}}a_{l_1}\cdots a_{l_{m-1}}\Big{)}_{1\le p,k\le n},
\end{eqnarray*}
where $f(z)=(f_1(z),f_2(z),\cdots,f_n(z))', a=(a_1,a_2,\cdots,a_n)'\in \mathbb{C}^n$.

Next, some mappings are defined as follows.

{\bf Defintion 1.1\cite{GK2003}}\quad Let $f:D\to \mathbb{C}$ be a holomorphic function. Then $f$ is convex if and only if $f'(0)\ne 0$ and
\begin{eqnarray*}
{\rm Re}\Big[1+\frac{zf''(z)}{f'(z)}\Big]\ge 0, z\in D.
\end{eqnarray*}

{\bf Defintion 1.2\cite{GK2003}}\quad Let $f:D\to \mathbb{C}$ be a holomorphic function. We say that $f$ is convex of order $\alpha$, $0\le\alpha<1$, if $f'(0)\ne 0$ and
\begin{eqnarray*}
{\rm Re}\Big[1+\frac{zf''(z)}{f'(z)}\Big]\ge\alpha, z\in D.
\end{eqnarray*}

{\bf Definition 1.3\cite{LL2010}}\quad Suppose that $\alpha \in [0,1)$, and $f:B\to X$ is a locally biholomorphic mapping. If
\begin{eqnarray*}
{\rm Re}\{T_x[(Df(x))^{-1}(D^{2}f(x)(x^2)+Df(x)x)]\}\ge{\alpha \Arrowvert x\Arrowvert}, x\in B,
\end{eqnarray*}
then $f$ is said to be a quasi-convex mapping of type B and order $\alpha$ on $B$. Take the family as $Q_{B}^{\alpha}(B)$. When $\alpha=0$, it is the family of quasi-convex mappings of type B.

{\bf Definition~1.4}\cite{LL2010}\quad Suppose that $\alpha \in [0,1)$, and $f:B\to X$ is a locally biholomorphic mapping. If
\begin{eqnarray*}
{\rm Re}\{T_x[(Df(x))^{-1}f(x)]\}\ge{\alpha \Arrowvert x\Arrowvert}, x\in B,
\end{eqnarray*}
then $f$ is said to be an almost starlike of order $\alpha$ on $B$.

{\bf Definition~1.5}\cite{LL2005}\quad Suppose $f\in H(B)$. We say that $x = 0$ is the zero of order $k$ of
$f(x)$ if $f(0) = 0,\cdots,D^{k-1}f(0) = 0$, but $D^kf(0)\ne 0$, where $k\in \mathbb{N}_+$. Note that the definition is
the same as that in the case $X =\mathbb{C}$.

{\bf Definition~1.6}\cite{IS1986,MK2012}\quad Suppose that $L:X^m\to\mathbb{C}$ is a continuous $m-$linear form, if
\begin{eqnarray*}
L(x_1,x_2,\cdots,x_m)=L(x_{\sigma(1)},x_{\sigma(2)},\cdots,x_{\sigma(m)}),\quad\forall x_1,\cdots,x_m\in X,
\end{eqnarray*}
for any $x_1,\cdots,x_n$ in $X$ and any permutation $\sigma$ of the first $m$ natural numbers, and
\begin{eqnarray*}
\|L\|=\sup{\{|L(x_1,x_2,\cdots,x_m)|:\|x_1\|\le1,\cdots,\|x_m\|\le1\}},
\end{eqnarray*}
then $L$ is said to be continuous symmetric $m-$linear form.
Denote $\mathcal{L}^s(^mX)$ to be the space of all continuous symmetric $m-$linear forms.

{\bf Definition~1.7}\cite{IS1986,MK2012}\quad Suppose that $L:X^m\to\mathbb{C}$ is a continuous symmetric $m-$linear form, if
\begin{eqnarray*}
P(x)=L(x,\cdots,x),\quad \forall x\in X,
\end{eqnarray*}
then $P:X\to \mathbb{C}$ is said to be continuous homogeneous polynomial of degree $m$. Let $\|P\|=\sup{\{|P(x)|:\|x\|\le1\}}$. Take the family as $\mathcal{P}^{s}(^mX)$.

For the sake of convenience, we let $\hat{L}=P$.

\setcounter{equation}{0}
\section{The inequalities of homogeneous expansion for biholomorphic quasi-convex mappings of type B and almost starlike mappings of order $\alpha$}

\qquad In order to prove the main results in this section, we need the following lemmas.

{\bf Lemma~2.1}\cite{GK2003}\quad If $f(z)=a_0+\sum_{n=1}^{\infty}a_{n}z^n\in H(D)$, and $f(D)\subset D$, then
\begin{eqnarray*}
|a_n|\le1-|a_0|^2, n=1,2,\cdots,
\end{eqnarray*}
when $n=1$, the above estimate is sharp.

{\bf Lemma~2.2}\quad  Let $p(z)=1+\sum_{n=1}^{\infty}b_{n}z^n\in H(D)$, and $\alpha \in [0,1)$. If ${\rm Re}\, p(z)\ge\alpha, z\in D$, then
\begin{eqnarray}
&&\Big{|}b_2-\frac{b_{1}^2}{2(1-\alpha)}\Big{|}\le 2(1-\alpha)-\frac{|b_{1}|^2}{2(1-\alpha)},\label{l21}\\
&&\Big{|}b_3-\frac{b_{1}b_2}{1-\alpha}+\frac{b_{1}^3}{4(1-\alpha)^2}\Big{|}\le 2(1-\alpha)-\frac{|b_{1}|^2}{2(1-\alpha)}.\label{l22}
\end{eqnarray}
Furthermore, we have
\begin{eqnarray}
|b_2|\le 2(1-\alpha),\label{l23}
\end{eqnarray}
and
\begin{eqnarray}
\Big{|}b_2-\frac{b_{1}^2}{1-\alpha}\Big{|}\le 2(1-\alpha).\label{l24}
\end{eqnarray}

{\bf Proof}\quad Let $h(z)=\frac{1-\frac{p(z)-\alpha}{1-\alpha}}{1+\frac{p(z)-\alpha}{1-\alpha}}=\frac{1-p(z)}{1-2\alpha+p(z)}, z\in D$,
then $h(0)=0,h(D)\subset D,h(z)\in H(D)$, that is, $h(z)$ is just a $Schwarz$ function, moreover,
\begin{eqnarray*}
h(z)&=&\frac{-\sum_{n=1}^{\infty}b_{n}z^n}{2(1-\alpha)+\sum_{n=1}^{\infty}b_{n}z^n}\\
&=&\frac{-b_1}{2(1-\alpha)}z-
\Big{[}\frac{b_2}{2(1-\alpha)}-\frac{b_{1}^2}{4(1-\alpha)^2}\Big{]}z^2\\
&&-\Big{[}\frac{b_3}{2(1-\alpha)}-\frac{b_{1}b_2}{2(1-\alpha)^2}+\frac{b_{1}^3}{8(1-\alpha)^3}\Big{]}z^3+\cdots.
\end{eqnarray*}

According to Schwarz Lemma, we have
\begin{eqnarray*}
|h(z)/z|=1 \  or \ |h(z)/z|<1.
\end{eqnarray*}

When $|h(z)/z|=1$, we have
\begin{eqnarray}
\Big{|}\frac{b_{1}}{2(1-\alpha)}\Big{|}=1 , \frac{b_2}{2(1-\alpha)}-\frac{b_{1}^2}{4(1-\alpha)^2}=0 ,
\frac{b_3}{2(1-\alpha)}-\frac{b_{1}b_2}{2(1-\alpha)^2}+\frac{b_{1}^3}{8(1-\alpha)^3}=0.\label{l25}
\end{eqnarray}

When $|h(z)/z|<1$, we have
\begin{eqnarray}
&&\Big{|}\frac{b_2}{2(1-\alpha)}-\frac{b_{1}^2}{4(1-\alpha)^2}\Big{|}\le1-\Big{|}\frac{b_1}{2(1-\alpha)}\Big{|}^2 ,\label{l26}\\
&&\Big{|}\frac{b_3}{2(1-\alpha)}-\frac{b_{1}b_2}{2(1-\alpha)^2}+\frac{b_{1}^3}{8(1-\alpha)^3}\Big{|}\le1-\Big{|}\frac{b_1}{2(1-\alpha)}\Big{|}^2.\label{l27}
\end{eqnarray}

By (\ref{l25}), (\ref{l26}) and (\ref{l27}), we obtain
\begin{eqnarray*}
&&\Big{|}b_2-\frac{b_{1}^2}{2(1-\alpha)}\Big{|}\le 2(1-\alpha)-\frac{|b_{1}|^2}{2(1-\alpha)},\\
&&\Big{|}b_3-\frac{b_{1}b_2}{1-\alpha}+\frac{b_{1}^3}{4(1-\alpha)^2}\Big{|}\le 2(1-\alpha)-\frac{|b_{1}|^2}{2(1-\alpha)}.
\end{eqnarray*}

Finally, we have
$$|b_2|\leq \Big{|}b_2-\frac{b_{1}^2}{2(1-\alpha)}\Big{|}+\frac{|b_{1}|^2}{2(1-\alpha)}\leq 2(1-\alpha),$$
and
\begin{eqnarray*}
\Big{|}b_2-\frac{b_{1}^2}{1-\alpha}\Big{|}\leq\Big{|}b_2-\frac{b_{1}^2}{2(1-\alpha)}\Big{|}+\frac{|b_{1}|^2}{2(1-\alpha)}\le 2(1-\alpha).
\end{eqnarray*}
This completes the proof.

{\bf Lemma~2.3}\cite{LL2009}\quad If $f(x):B\to X$ is a normalized locally biholomorphic mapping, and $g(x)=(Df(x))^{-1}\Big(D^2f(x)(x^2)+Df(x)x\Big)$, then
\begin{eqnarray*}
&&\frac{D^2g(0)(x^2)}{2!}=2\frac{D^2f(0)(x^2)}{2!},\\
&&\frac{D^3g(0)(x^3)}{3!}=6\frac{D^3f(0)(x^3)}{3!}-4\frac{D^2f(0)}{2!}\Big{(}x,\frac{D^2f(0)(x^2)}{2!}\Big{)}.
\end{eqnarray*}

{\bf Lemma~2.4}\cite{LL2005}\quad If $f(x):B\to X$ is a normalized locally biholomorphic mapping, and $g(x)=(Df(x))^{-1}f(x)$, then
\begin{eqnarray*}
\frac{D^2g(0)(x^2)}{2!}&=&-\frac{D^2f(0)(x^2)}{2!},\\
\frac{D^3g(0)(x^3)}{3!}&=&-2\frac{D^3f(0)(x^3)}{3!}+2\frac{D^2f(0)}{2!}\Big{(}x,\frac{D^2f(0)(x^2)}{2!}\Big{)},\\
\frac{D^4g(0)(x^4)}{4!}&=&-3\frac{D^4f(0)(x^4)}{4!}+3\frac{D^3f(0)}{3!}\Big{(}x^2,\frac{D^2f(0)(x^2)}{2!}\Big{)}\\
&&+4\frac{D^2f(0)}{2!}\Big{(}x,\frac{D^3f(0)(x^3)}{3!}\Big{)}-4\frac{D^2f(0)}{2!}\Big{(}x,\frac{D^2f(0)}{2!}(x,\frac{D^2f(0)(x^2)}{2!})\Big{)}.
\end{eqnarray*}

Now we establish the sharp inequalities of homogeneous expansion for biholomorphic quasi-convex mappings of type B and order $\alpha$, and the almost starlike mappings of order $\alpha$.

{\bf Theorem~2.1}\quad If $f(x)\in Q_{B}^{\alpha}(B),\ 0\le \alpha<1$, then
\begin{eqnarray*}
\Big{|}T_x\Big{(}\frac{D^3f(0)(x^3)}{3!}\Big{)}-\frac{2}{3}T_x\Big{(}\frac{D^2f(0)}{2!}(x,\frac{D^2f(0)(x^2)}{2!})\Big{)}\Big{|}\le \frac{1-\alpha}{3}\|x\|^3, x\in B.
\end{eqnarray*}
The above inequality is sharp.

{\bf Proof}\quad For fixed $x\in B\backslash\{0\}$, let $x_0=\frac{x}{\|x\|}$. Define $p(\xi)=\frac{T_{x_0}(g(\xi x_0))}{\xi}, \xi\in D$,  where
\begin{eqnarray*}
g(x)=(Df(x))^{-1}(D^2f(x)(x^2)+Df(x)x).
\end{eqnarray*}

By the hypothesis of Theorem~2.1, we have $p(\xi)\in H(D)$, ${\rm Re} p(\xi)\ge \alpha$, $p(0)=1$, and
\begin{eqnarray*}
p(\xi)=1+\frac{T_{x_0}(D^2g(0)(x_{0}^2))}{2!}\xi+\frac{T_{x_0}(D^3g(0)(x_{0}^3))}{3!}\xi^2+\cdots.
\end{eqnarray*}

According to (\ref{l23}) of Lemma~2.2, we obtain
\begin{eqnarray*}
\Big{|}T_{x_0}(\frac{D^3g(0)(x_{0}^3)}{3!})\Big{|}\le 2(1-\alpha),\quad x_0\in \partial B.
\end{eqnarray*}

From Lemma~2.3, we conclude that
\begin{eqnarray*}
\frac{D^3g(0)(x_{0}^3)}{3!}&=&6\frac{D^3f(0)(x_{0}^3)}{3!}-4\frac{D^2f(0)}{2!}\Big{(}x_0,\frac{D^2f(0)(x_{0}^2)}{2!}\Big{)}\\
&=&6\Big{[}\frac{D^3f(0)(x_{0}^3)}{3!}-\frac{2}{3}\frac{D^2f(0)}{2!}\Big{(}x_0,\frac{D^2f(0)(x_{0}^2)}{2!}\Big{)}\Big{]}.
\end{eqnarray*}

That is
\begin{eqnarray*}
\Big{|}T_{x_0}\Big{(}\frac{D^3f(0)(x_{0}^3)}{3!}\Big{)}-\frac{2}{3}T_{x_0}\Big{(}\frac{D^2f(0)}{2!}(x_0,\frac{D^2f(0)(x_{0}^2)}{2!})\Big{)}\Big{|}
\le \frac{1-\alpha}{3}.
\end{eqnarray*}

Hence,
\begin{eqnarray*}
\Big{|}T_x\Big{(}\frac{D^3f(0)(x^3)}{3!}\Big{)}-\frac{2}{3}T_x\Big{(}\frac{D^2f(0)}{2!}(x,\frac{D^2f(0)(x^2)}{2!})\Big{)}\Big{|}\le \frac{1-\alpha}{3}\|x\|^3\, \mbox{ for }\,  x\in B.
\end{eqnarray*}

Finally, it is easy to check that the function
\begin{eqnarray*}
f(x)=\left\{
\begin{array}{lll}
\frac{1-(1-T_u(x))^{2\alpha-1}}{2\alpha-1}x, & \alpha\in[0,1)\backslash\{\frac{1}{2}\},\\
-x\log{(1-T_u(x))}, & \alpha=\frac{1}{2}
\end{array}
\right.
\end{eqnarray*}
satisfies the condition of Theorem~2.1, where $x\in B, u\in\partial B$. We set $x=ru, \|u\|=1, 0\le r<1$. By a direct computation, we obtain that
\begin{eqnarray*}
\Big{|}T_x\Big{(}\frac{D^3f(0)(x^3)}{3!}\Big{)}-\frac{2}{3}T_x\Big{(}\frac{D^2f(0)}{2!}(x,\frac{D^2f(0)(x^2)}{2!})\Big{)}\Big{|}
= \frac{1-\alpha}{3}r^3.
\end{eqnarray*}
Hence, the inequality in the Theorem~2.1 is sharp.

{\bf Remark~2.1}\quad Setting $X=\mathbb{C}^n, B=D^n,\alpha=0$ in Theorem~2.1, then we can deduce Theorem~2.1 in \cite{LL2009}.

Applying the similarly method as in the proof of Theorem~2.1, we can prove the following theorem.

{\bf Theorem~2.2}\quad If $f(x)\in Q_{B}^{\alpha}(D^n),0\le \alpha<1$, then
\begin{eqnarray*}
\Big{\|}\frac{D^3f(0)(x^3)}{3!}-\frac{2}{3}\frac{D^2f(0)}{2!}(x,\frac{D^2f(0)(x^2)}{2!})\Big{\|}\le \frac{1-\alpha}{3}\|x\|^3, \forall x\in D^n.
\end{eqnarray*}
The above inequality is sharp.

Setting $n=1$ in Theorem~2.2, we obtain Corollary~2.2.

{\bf Corollary~2.2}\quad If $f(z)=z+\sum_{n=2}^{\infty}a_{n}z^n$ is a univalent convex function of order $\alpha$ on the unit disk $D$ in the $\mathbb{C}$, $0\le\alpha<1$, then
\begin{eqnarray*}
|a_3-\frac{2}{3}a_{2}^2|\le\frac{1-\alpha}{3}.
\end{eqnarray*}
The above inequality is sharp.

Setting $\alpha=0$ in Corollary~2.2, we obtain Corollary~E.

{\bf Theorem~2.3}\quad If $f(x)$ is an almost starlike of order $\alpha$ on $B$, $0\le \alpha<1$, then
\begin{eqnarray*}
&&\Big{|}2T_x\Big{(}\frac{D^3f(0)(x^3)}{3!}\Big{)}\|x\|-2T_x\Big{(}\frac{D^2f(0)}{2!}(x,\frac{D^2f(0)(x^2)}{2!})\Big{)}\|x\|
+\frac{[T_x(\frac{D^2f(0)(x^2)}{2!})]^2}{1-\alpha}\Big{|}\\
&\leq& 2(1-\alpha)\|x\|^4, \forall x\in B.
\end{eqnarray*}
The above inequality is sharp.

{\bf Proof}\quad For fixed $x\in B\backslash\{0\}$, let $x_0=\frac{x}{\|x\|}$. Define $p(\xi)=\frac{T_{x_0}(g(\xi x_0))}{\xi}, \xi\in D$,  where  $$g(x)=(Df(x))^{-1}f(x).$$
By the hypothesis of Theorem~2.1, we have $p(\xi)\in H(D)$, ${\rm Re}p(\xi)=Re\frac{1}{|\xi|}T_{\xi x_0}(g(\xi x_0))\ge \alpha$, $$p(0)=\lim_{\xi\to0}\frac{1}{\xi}T_{x_0}(g(\xi x_0))=T_{x_0}(x_0)=1,$$
and
\begin{eqnarray*}
p(\xi)=1+\frac{T_{x_0}(D^2g(0)(x_{0}^2))}{2!}\xi+\frac{T_{x_0}(D^3g(0)(x_{0}^3))}{3!}\xi^2+\cdots.
\end{eqnarray*}

Consequently, by (\ref{l24}) in Lemma~2.2, we obtain
\begin{eqnarray*}
\Big{|}T_{x_0}\Big{(}\frac{D^3g(0)(x_{0}^3)}{3!}\Big{)}-\frac{1}{1-\alpha}\Big{[}T_{x_0}\Big{(}\frac{D^2g(0)(x_{0}^2)}{2!}\Big{)}\Big{]}^2\Big{|}\le 2(1-\alpha).
\end{eqnarray*}

From Lemma~2.4, we conclude that
\begin{eqnarray*}
&&\frac{D^2g(0)(x_{0}^2)}{2!}=-\frac{D^2f(0)(x_{0}^2)}{2!},\\
&&\frac{D^3g(0)(x_{0}^3)}{3!}=-2\frac{D^3f(0)(x_{0}^3)}{3!}+2\frac{D^2f(0)}{2!}\Big{(}x_0,\frac{D^2f(0)(x_{0}^2)}{2!}\Big{)},
\end{eqnarray*}
that is
\begin{eqnarray*}
\Big{|}2T_{x_0}\Big{(}\frac{D^3f(0)(x_{0}^3)}{3!}\Big{)}-2T_{x_0}\Big{(}\frac{D^2f(0)}{2!}(x_0,\frac{D^2f(0)(x_{0}^2)}{2!})\Big{)}
+\frac{[T_{x_0}(\frac{D^2f(0)(x_{0}^2)}{2!})]^2}{1-\alpha}\Big{|}\le 2(1-\alpha), \forall x_0\in \partial B,
\end{eqnarray*}

For $x_0=\frac{x}{\|x\|}$, so we have $T_x(\cdot)=T_{x_0}(\cdot)$. Thus we obtain
\begin{eqnarray*}
\Big{|}2T_x\Big{(}\frac{D^3f(0)}{3!}(\frac{x^3}{\|x\|^3})\Big{)}-2T_x\Big{(}\frac{D^2f(0)}{2!}(\frac{x}{\|x\|},
\frac{D^2f(0)}{2!}(\frac{x^2}{\|x\|^2}))\Big{)}+\frac{[T_x(\frac{D^2f(0)}{2!}(\frac{x^2}{\|x\|^2}))]^2}{1-\alpha}\Big{|}\le 2(1-\alpha), \forall x\in B,
\end{eqnarray*}
that is,
\begin{eqnarray*}
\Big{|}2T_x\Big{(}\frac{D^3f(0)(x^3)}{3!}\Big{)}\|x\|-2T_x\Big{(}\frac{D^2f(0)}{2!}(x,\frac{D^2f(0)(x^2)}{2!})\Big{)}\|x\|
+\frac{[T_x(\frac{D^2f(0)(x^2)}{2!})]^2}{1-\alpha}\Big{|}\le 2(1-\alpha)\|x\|^4, \forall x\in B.
\end{eqnarray*}

Finally, it is not difficult to check that the function
\begin{eqnarray*}
f(x)=\left\{
\begin{array}{lll}
\frac{x}{(1-(1-2\alpha)T_u(x))^{\frac{2(1-\alpha)}{1-2\alpha}}},& \alpha\in[0,1)\backslash\{\frac{1}{2}\},\\
xe^{T_u(x)}, & \alpha=\frac{1}{2}
\end{array}
\right.
\end{eqnarray*}
satisfies the hypothesis of Theorem~2.3, where $u\in\partial B$. We set $x=ru, \|u\|=1, 0\le r<1$. By direct computation, we obtain that
\begin{eqnarray*}
\Big{|}2T_x\Big{(}\frac{D^3f(0)(x^3)}{3!}\Big{)}\|x\|-2T_x\Big{(}\frac{D^2f(0)}{2!}(x,\frac{D^2f(0)(x^2)}{2!})\Big{)}\|x\|+
\frac{[T_x(\frac{D^2f(0)(x^2)}{2!})]^2}{1-\alpha}\Big{|}= 2(1-\alpha)\|r\|^4.
\end{eqnarray*}
Hence, the inequality in Theorem~2.3 is sharp.

Setting $\alpha=0$ in Theorem~2.3, we can obtain the following corollary.

{\bf Corollary~2.3}\cite{LL2013}\quad If $f(x)\in S^*(B)$. Then
\begin{eqnarray*}
\Big{|}2T_x\Big{(}\frac{D^3f(0)(x^3)}{3!}\Big{)}\|x\|-2T_x\Big{(}\frac{D^2f(0)}{2!}(x,\frac{D^2f(0)(x^2)}{2!})\Big{)}\|x\|
+[T_x(\frac{D^2f(0)(x^2)}{2!})]^2\Big{|}\le 2\|x\|^4, \forall x\in B,
\end{eqnarray*}
The above inequality is sharp.

When $\alpha=0, \ n=1$, we obtain Corollary~F by Corollary~2.3.

\setcounter{equation}{0}
\section{Applications of inequalities}
\quad In this section, with the sharp inequalities in Section 2, we will establish the sharp estimates of the third and fourth homogeneous expansions for the normalized locally biholomorphic quasi-convex mappings of type B and order $\alpha$ and almost starlike mappings of order $\alpha$ on $D^n$ in $\mathbb{C}^n$.

In order to establish the main results in this section, we need the following lemmas.

{\bf Lemma~3.1}\cite{LL2009}\quad Suppose $0\le\alpha<1, f:D^n\to X$ is a normalized locally biholomorphic mapping, then $f\in Q_{B}^{\alpha}(D^n)$ if and only if
\begin{eqnarray*}
\mbox{Re}\Big{\{}\frac{g_j(z)}{z_j}\Big{\}}\ge \alpha, \forall z\in D^n,
\end{eqnarray*}
where $g(z)=(g_1(z),g_2(z),\cdots,g_n(z))'=(Df(z))^{-1}(D^2f(z)(z^2)+Df(z)z)$ is a column vector in $\mathbb{C}^n$, $|z_j|=\|z\|=\max_{1\le k\le n}\{|z_k|\}$.

{\bf Lemma~3.2}\cite{LL2009}\quad Suppose $0\le\alpha<1, f:D^n\to X$ is a normalized locally biholomorphic mapping, then $f$ is an almost starlike mapping of order $\alpha$ if and only if
\begin{eqnarray*}
\mbox{Re}\Big{\{}\frac{g_j(z)}{z_j}\Big{\}}\ge \alpha, \forall z\in D^n,
\end{eqnarray*}
where $g(z)=(g_1(z),g_2(z),\cdots,g_n(z))'=(Df(z))^{-1}f(z))$ is a column vector in $\mathbb{C}^n$, $|z_j|=\|z\|=\max_{1\le k\le n}\{|z_k|\}$.

{\bf Lemma~3.3}\cite{LL2005}\quad Suppose $g\in H(D^n), g(0)=0, Dg(0)=I,\alpha\in[0,1)$. If $\mbox{Re}\frac{g_j(z)}{z_j}\ge \alpha, z\in D^n$, where $|z_j|=\|z\|=\max_{1\le k\le n}\{|z_k|\}$, then
\begin{eqnarray*}
\frac{\|D^mg(0)(z^m)\|}{m!}\le2(1-\alpha)\|z\|^m, z\in D^n, m=2,3,\cdots.
\end{eqnarray*}
The above estimate is sharp.

{\bf Lemma~3.4}\quad If
\begin{eqnarray*}
&&\Big\|(z_{1}^q(\sum_{l=1}^{n}a_{p1l}z_l^p),z_{2}^q(\sum_{l=1}^{n}a_{p2l}z_l^p),\cdots,z_{n}^q(\sum_{l=1}^{n}a_{pnl}z_l^p))'\Big\|\le C_0\|z\|^m,\\
&&z\in D^n,p+q=m,p,q\in\mathbb{N}^+,m=2,3,\cdots,
\end{eqnarray*}
where each $a_{pkl}(p=1,2,\cdots,m-1,k,l=1,2,\cdots,n)$ is a complex number independent of $z_k(k=1,2,\cdots,n)$, $\|z\|=\max_{1\le k\le n}\{|z_k|\}$, $C_0$ is a non-negative real constant. Then
\begin{eqnarray*}
M=\max_{1\le k\le n}\{\sum_{l=1}^n|a_{pkl}|\}\le C_0,p=1,2,\cdots,m-1.
\end{eqnarray*}

{\bf Proof}\quad For every $z\in D^n\backslash\{0\}$, taking into the hypothesis of Lemma~3.3, we have
\begin{eqnarray*}
\Big{|}\frac{z_{k}^q}{\|z\|^q}\Big{(}\sum_{l=1}^{n}{a_{pkl}\big{(}\frac{z_l}{\|z\|}\big{)}^p}\Big{)}\Big{|}\le C_0.
\end{eqnarray*}
Especially, for each $k$, if $a_{pkl}\ne0$, taking $z_l=e^{-i\frac{\arg{a_{pkl}}}{p}}\|z\|, l=1,2,\cdots,n$, where $i^2=-1$. Then, we conclude that
\begin{eqnarray*}
\sum_{l=1}^{n}{|a_{pkl}|}\le C_0, k=1,2,\cdots,n,
\end{eqnarray*}
that is
\begin{eqnarray*}
M=\max_{1\le k\le n}{\sum_{l=1}^n{\{|a_{pkl}|}\}}\le C_0,p=1,2,\cdots,m-1.
\end{eqnarray*}
This completes this proof.

{\bf Lemma~3.5}\cite{IS1986}\quad Suppose $L\in\mathcal{L}^s(^m \ l_{n}^{\infty})$, $\hat{L}\in\mathcal{P}^s(^m\ l_{n}^{\infty})$, where $\mathcal{L}^s(^m \ l_{n}^{\infty}),\mathcal{P}^s(^m\ l_{n}^{\infty})$ are defined by Definitions~1.6 and 1.7, then
\begin{eqnarray*}
&&\|L\|=\|\hat{L}\|,n=2,\\
&&\|L\|\le m^{m/2}(m+1)^{(m+1)/2}/2^{m}m!\|\hat{L}\|,n\ge3,n\in\mathbb{N}_+.
\end{eqnarray*}

Especially, when $m=2$, according to Lemma~3.5, we have

{\bf Lemma~3.6}\quad Suppose $L\in\mathcal{L}^s(^2 \ l_{n}^{\infty})$, $\hat{L}\in\mathcal{P}^s(^2\ l_{n}^{\infty})$, where $\mathcal{L}^s(^2 \ l_{n}^{\infty}),\mathcal{P}^s(^2\ l_{n}^{\infty})$ are defined by Definitions~1.6 and 1.7, then
\begin{eqnarray*}
&&\|L\|=\|\hat{L}\|,n=2,\\
&&\|L\|\le \frac{3}{4}\sqrt{3}\|\hat{L}\|,n\ge3,n\in\mathbb{N}_+.
\end{eqnarray*}

By Lemma~3.6, we can obtain the following lemma.

{\bf Lemma~3.7}\quad Suppose $f:D^n\to X$ is a holomorphic mapping. Define
\begin{eqnarray*}
&&L(x,y)=D^2f(0)(x,y),x,y\in X,\ \ \text{ÇÒ}\|L\|=\sup{\{\|L(x,y)\|:\|x\|\le1,\|y\|\le1\}},\\
&&\hat{L}(x,x)=D^2f(0)(x^2),x\in X,\ \ \ \ \text{ÇÒ}\|\hat{L}\|=\sup{\{\|\hat{L}(x,x)\|:\|x\|\le1\}},
\end{eqnarray*}
where $\|x\|=\max_{1\le k\le n}{\{|x_k|\}},\  \|y\|=\max_{1\le k\le n}{\{|y_k|\}}$. If $L$ and $\hat{L}$ are bounded linear operators, then
\begin{eqnarray*}
&&\|L\|=\|\hat{L}\|,\, \mbox{ for }\, n=2,\\
&&\|L\|\le\frac{3}{4}\sqrt{3}\|\hat{L}\|,\, \mbox{ for }\, n\ge3,\, n\in\mathbb{N}_+.
\end{eqnarray*}
Furthermore
\begin{eqnarray*}
\|D^2f(0)(x,y)\|\le\|L\|\|x\|\|y\|,\, \mbox{ for }\, x,y\in X.
\end{eqnarray*}

{\bf Proof}\quad Since $L$ and $\hat{L}$ are bounded linear operators, according to Lemma~3.6, we have
\begin{eqnarray*}
&&\|L\|=\sup{\{\|L(x,y)\|:\|x\|\le1,\|y\|\le1\}}=\sup_{x,y\ne0,x,y\in X}{\frac{\|D^2f(0)(x,y)\|}{\|x\|\|y\|}},\\
&&\|\hat{L}\|=\sup{\{\|\hat{L}(x,x)\|:\|x\|\le1\}}=\sup_{x\ne0,x\in X}{\frac{\|D^2f(0)(x,x)\|}{\|x\|^2}},
\end{eqnarray*}
then
\begin{eqnarray*}
\|D^2f(0)(x,y)\|\le\|L\|\|x\|\|y\|,x,y\in X,\\
\|D^2f(0)(x,x)\|\le\|\hat{L}\|\|x\|^2,x,\in X.
\end{eqnarray*}
Hence, when $\|x\|\le1,\|y\|\le1$, $\|D^2f(0)(x,y)\|,\|D^2f(0)(x,x)\|$ are all bounded, where
\begin{eqnarray*}
&&\|D^2f(0)(x,y)\|=\max_{1\le k\le n}{\{|D^2f_{k}(0)(x,y)|,\|x\|\le1,\|y\|\le1\}},\\
&&\|D^2f(0)(x,x)\|=\max_{1\le k\le n}{\{|D^2f_{k}(0)(x^2)|,\|x\|\leq 1\}}.
\end{eqnarray*}
From Lemma~3.6, we have
\begin{eqnarray*}
&&L(x,y)=\Big{(}\sum_{i,j=1}^{n}{a_{ij}^{1}x_{i}y_{j}},\cdots,\sum_{i,j=1}^{n}{a_{ij}^{n}x_{i}y_{j}}\Big{)}',\\
&&\hat{L}(x,x)=\Big{(}\sum_{i,j=1}^{n}{a_{ij}^{1}x_{i}x_{j}},\cdots,\sum_{i,j=1}^{n}{a_{ij}^{n}x_{i}x_{j}}\Big{)}',
\end{eqnarray*}
where $a_{ij}^{k}=\frac{\partial^2f_{k}(0)}{\partial x_i\partial x_j},\ i,j,k=1,2,\cdots,n$. Let
\begin{eqnarray*}
L(x,y)=(L_1(x,y),\cdots,L_{n}(x,y))',\\
\hat{L}(x,x)=(\hat{L}_{1}(x,x),\cdots,\hat{L_n}(x,x))',
\end{eqnarray*}
where $L_k=\sum_{i,j=1}^{n}{a_{ij}^{k}x_{i}y_{j}},\ \hat{L}_{k}=\sum_{i,j=1}^{n}{a_{ij}^{k}x_{i}x_{j}},\ k=1,2,\cdots,n$. By Definitions~1.6 and 1.7, it is easy to obtain that $L_k\in\mathcal{L}^s(^2\  l_{n}^{\infty})$, $\hat{L}_{k}\in\mathcal{P}^s(^2\ l_{n}^{\infty}),\ k=1,2,\cdots,n$. then $L_k$ and $\hat{L}_{k}$  satisfy the conditions of Lemma~3.6. Therefore,
\begin{eqnarray}
\|L_k\|&=&\|\hat{L}_{k}\|,\, \mbox{ for }\, n=2,\label{liu31}\\
\|L_k\|&\le&\frac{3}{4}\sqrt{3}\|\hat{L}_{k}\|,\, \mbox{ for }\, n\ge3,\, n\in\mathbb{N}_+.\label{liu32}
\end{eqnarray}
Because
\begin{eqnarray*}
\|L\|&=&\sup_{\|x\|\le1,\|y\|\le1}{\{\|L(x,y)\|\}}\\
&=&\sup_{\|x\|\le1,\|y\|\le1}{\max_{1\le k\le n}{\{|L_{k}(x,y)|\}}}\\
&\ge&\sup_{\|x\|\le1,\|y\|\le1}{\{|L_{k}(x,y)|\}}\\
&=&\|L_{k}\|,k=1,2,\cdots,n,
\end{eqnarray*}
we have
\begin{eqnarray*}
\|L\|\ge\max_{1\le k\le n}{\{\|L_{k}\|\}}.
\end{eqnarray*}
On the other hand,
\begin{eqnarray*}
\|L\|&=&\sup_{\|x\|\le1,\|y\|\le1}{\{\|L(x,y)\|\}}\\
&=&\sup_{\|x\|\le1,\|y\|\le1}{\max_{1\le k\le n}{\{|L_{k}(x,y)|\}}}\\
&\le&\sup_{\|x\|\le1,\|y\|\le1}{\max_{1\le k\le n}{\sup_{\|x\|\le1,\|y\|\le1}{\{|L_{k}(x,y)|\}}}}\\
&=&\sup_{\|x\|\le1,\|y\|\le1}{\max_{1\le k\le n}{\{\|L_{k}\|\}}}\\
&=&\max_{1\le k\le n}{\{\|L_{k}\|\}}.
\end{eqnarray*}
Hence,
\begin{eqnarray*}
\|L\|=\max_{1\le k\le n}{\{\|L_{k}\|\}}.
\end{eqnarray*}
Similarly, we have
\begin{eqnarray*}
\|\hat{L}\|=\max_{1\le k\le n}{\{\|\hat{L}_{k}\|\}}.
\end{eqnarray*}
According to (\ref{liu31}) and (\ref{liu32}), we obtain
\begin{eqnarray*}
&&\|L\|=\|\hat{L}\|,n=2,\\
&&\|L\|\le\frac{3}{4}\sqrt{3}\|\hat{L}\|,n\ge3,n\in\mathbb{N}_+.
\end{eqnarray*}
This completes the proof.

{\bf Lemma~3.8}\quad Assume that $\alpha\in[0,\frac{37-\sqrt{505}}{72}]$, and $$h(x)=\frac{12\alpha^2-10\alpha+1}{4(1-\alpha)^2}x^3-\frac{1}{2(1-\alpha)}x^2+(5-7\alpha)x+2(1-\alpha),$$
then $h(x)$ is strictly increasing on $x\in [0,2(1-\alpha)]$, and
$$\max_{x\in[0,2(1-\alpha)]}{\{h(x)\}}=h(2(1-\alpha))=(1-\alpha)(3-4\alpha)(4-6\alpha).$$

{\bf Proof}\quad Since $h(x)=\frac{12\alpha^2-10\alpha+1}{4(1-\alpha)^2}x^3-\frac{1}{2(1-\alpha)}x^2+(5-7\alpha)x+2(1-\alpha)$, we have
\begin{eqnarray*}
h'(x)&=&\frac{3(12\alpha^2-10\alpha+1)}{4(1-\alpha)^2}x^2-\frac{1}{1-\alpha}x+5-7\alpha\\
&=&\frac{3(12\alpha^2-10\alpha+1)}{4(1-\alpha)^2}\Big{(}x-\frac{2(1-\alpha)}{3(12\alpha^2-10\alpha+1)}\Big{)}^2\\
&&+5-7\alpha-\frac{1}{3(12\alpha^2-10\alpha+1)}.
\end{eqnarray*}

When $\alpha\in[0,\frac{15-\sqrt{154}}{36})$,
\begin{eqnarray*}
12\alpha^2-10\alpha+1>0,\ \ 0<\frac{2(1-\alpha)}{3(12\alpha^2-10\alpha+1)}<2(1-\alpha),\\
 h'(\frac{2(1-\alpha)}{3(12\alpha^2-10\alpha+1)})=5-7\alpha-\frac{1}{3(12\alpha^2-10\alpha+1)}>0,
\end{eqnarray*}
then $h'(x)\geq  h'(\frac{2(1-\alpha)}{3(12\alpha^2-10\alpha+1)})>0$ for each $x\in[0,2(1-\alpha)]$. Hence $h(x)$ is strictly increasing on $x\in [0,2(1-\alpha)]$.\\

When $\alpha\in[\frac{15-\sqrt{154}}{36},\frac{10-\sqrt{52}}{24})$,
\begin{eqnarray*}
12\alpha^2-10\alpha+1>0,\ \ \frac{2(1-\alpha)}{3(12\alpha^2-10\alpha+1)}\ge2(1-\alpha),\ \ h'(0)=5-7\alpha>0,\ \ h'(2(1-\alpha))>0,
\end{eqnarray*}
then $h'(x)>0$, for every $x\in[0,2(1-\alpha)]$. So $h(x)$ is strictly increasing on $x\in [0,2(1-\alpha)]$.

When $\alpha=\frac{10-\sqrt{52}}{24}$,
\begin{eqnarray*}
h'(x)=-\frac{1}{1-\alpha}x+5-7\alpha, h'(2(1-\alpha))=3-7\alpha>0,
\end{eqnarray*}
then $h'(x)>0$, for $\forall x\in[0,2(1-\alpha)]$, $h(x)$ is strictly increasing on $x\in [0,2(1-\alpha)]$.

When $\alpha\in(\frac{10-\sqrt{52}}{24},\frac{37-\sqrt{505}}{72}]$,
\begin{eqnarray*}
12\alpha^2-10\alpha+1<0,\ \ \frac{2(1-\alpha)}{3(12\alpha^2-10\alpha+1)}<0,\ \ h'(2(1-\alpha))\ge0,
\end{eqnarray*}
then $h'(x)>0$, for $\forall x\in[0,2(1-\alpha)]$, $h(x)$ is strictly increasing on $x\in [0,2(1-\alpha)]$. Hence,
for $\forall\alpha\in[0,\frac{37-\sqrt{505}}{36}]$, where $\frac{37-\sqrt{505}}{36}\approx0.2018$, $h(x)$ is strictly increasing on $x\in [0,2(1-\alpha)]$, and
\begin{eqnarray*}
\max_{x\in[0,2(1-\alpha)]}{\{h(x)\}}=h(2(1-\alpha))=(1-\alpha)(3-4\alpha)(4-6\alpha).
\end{eqnarray*}
This completes this proof.

{\bf Theorem~3.1}\quad  If $0\le\alpha<1$ and $f\in Q_{B}^{\alpha}(D^n)$, then for any $z\in D^n$, we have
\begin{eqnarray}
&&\frac{\Arrowvert D^3f(0)(z^3)\Arrowvert}{3!}\le \frac{(1-\alpha)(3-2\alpha)}{3}\Arrowvert z\Arrowvert^3,\, \mbox{ for }\,  n=2,\label{liu33}\\
&&\frac{\Arrowvert D^3f(0)(z^3)\Arrowvert}{3!}\le \frac{(1-\alpha)(\frac{3}{2}\sqrt{3}+1-\frac{3}{2}\sqrt{3}\alpha)}{3}\Arrowvert z\Arrowvert^3,\, \mbox{ for }\,  n\ge3,n\in\mathbb{N}_+.\label{liu34}
\end{eqnarray}
When $n=2$, the estimate (\ref{liu33}) is sharp.

{\bf Proof}\quad For any $z\in D^n\setminus\{0\}$, define $z_0=\frac{z}{\Arrowvert z\Arrowvert}$. By Lemmas~2.3 and 3.1, we have
\begin{eqnarray}
\Arrowvert D^2f(0)(z^{2}_0)\Arrowvert\le 2(1-\alpha),\|z_0\|=1,
\end{eqnarray}

Apparently, $\hat{L}(z^2)=D^2f(0)(z^{2})$ is a bounded linear operator. Then we have
\begin{eqnarray}
\|\hat{L}\|=\sup{\{\|\hat{L}(z^2)\|:\|z\|\le1\}}=\sup{\{\|\hat{L}(z^2)\|:\|z\|=1\}}.
\end{eqnarray}

According to $(3.3),\ (3.4)$, we obtain
\begin{eqnarray*}
\|\hat{L}\|=2(1-\alpha).
\end{eqnarray*}

Let $L(x,y)=D^2f(0)(x,y),x,y\in X$, by Lemma~3.6, we have
\begin{eqnarray*}
\|L\|&=&2(1-\alpha),n=2,\\
\|L\|&\le&\frac{3}{2}\sqrt{3}(1-\alpha),n\ge3,n\in\mathbb{N}_+£¬
\end{eqnarray*}

Let $w=\frac{D^2f(0)(z^2)}{2!}$, again by Lemma~3.7, we obtain
\begin{eqnarray*}
&&\|D^2f(0)(z,w)\|\le\|L\|\|z\|\|w\|\le2(1-\alpha)\|z\|\|w\|\le2(1-\alpha)^2\|z\|^2, z\in D^n,n=2,\\
&&\|D^2f(0)(z,w)\|\le\|L\|\|z\|\|w\|\le\frac{3}{2}\sqrt{3}(1-\alpha)\|z\|\|w\|\le\frac{3}{2}\sqrt{3}(1-\alpha)^2\|z\|^2, z\in D^n,n\ge3,n\in\mathbb{N}_+.
\end{eqnarray*}

According to Corollary~2.2, we conclude that when $n=2$,
\begin{eqnarray*}
\frac{\Arrowvert D^3f(0)(z^3)\Arrowvert}{3!}&\le&\frac{1-\alpha}{3}\Arrowvert z\Arrowvert^3+\frac{1}{3}\Arrowvert D^2f(0)(z,w)\Arrowvert\\
&\le&\frac{1-\alpha}{3}\Arrowvert z\Arrowvert^3+\frac{ 2(1-\alpha)^2}{3}\Arrowvert z\Arrowvert^3\\
&=&\frac{(1-\alpha)(3-2\alpha)}{3}\Arrowvert z\Arrowvert^3, z\in D^2.
\end{eqnarray*}

Similarly, when $n\ge3,n\in\mathbb{N}_+$,
\begin{eqnarray*}
\frac{\Arrowvert D^3f(0)(z^3)\Arrowvert}{3!}\le\frac{(1-\alpha)(\frac{3}{2}\sqrt{3}+1-\frac{3}{2}\sqrt{3}\alpha)}{3}\Arrowvert z\Arrowvert^3, z\in D^n.
\end{eqnarray*}

When $n=2$, it is easy to prove the function
\begin{eqnarray*}
f(z)=\left\{
\begin{array}{lll}
\Big(\frac{1-(1-z_1)^{2\alpha-1}}{2\alpha-1},\frac{1-(1-z_2)^{2\alpha-1}}{2\alpha-1})', & \alpha\in[0,1)\backslash\{\frac{1}{2}\},\\
(-\log(1-z_1),-\log(1-z_2))', & \alpha=\frac{1}{2}
\end{array}
\right.
\end{eqnarray*}
meets the conditions of Theorem 3.1. Taking $z=(r,0)'(0\le r<1)$, then
\begin{eqnarray*}
\frac{\Arrowvert D^3f(0)(z^3)\Arrowvert}{3!}=\frac{(1-\alpha)(3-2\alpha)}{3}r^3.
\end{eqnarray*}
This completes the proof.

{\bf Theorem~3.2}\quad Suppose $0\le\alpha<1$, $f\in Q_{B}^{\alpha}(D^n)$, and $\frac{D^2f_k(0)(z^2)}{2!}=z_k(\sum_{l=1}^na_{kl}z_l), k=1,2,\cdots$, where $a_{ml}=\frac{1}{2!}\frac{\partial^2f_m(0)}{\partial z_m\partial z_l}, m,l=1,2,\cdots,n$, then
\begin{eqnarray*}
\frac{\|D^3f(0)(z^3)\|}{3!}\le \frac{(1-\alpha)(3-2\alpha)}{3}\|z\|^3, z\in D^n.
\end{eqnarray*}
The above estimate is sharp.

{\bf Proof.}\quad For fixed $z\in D^n\setminus\{0\}$, define $z_0=\frac{z}{\Arrowvert z\Arrowvert}$. Through Lemmas~2.3, ~3.1 and ~3.2, we have
\begin{eqnarray*}
\Arrowvert D^2f(0)(z^{2}_0)\Arrowvert\le 2(1-\alpha).
\end{eqnarray*}

According to Lemma~3.4, we obtain
\begin{eqnarray}
\sum_{l=1}^{n}|a^{(j)}_{jl}|\le 2(1-\alpha).
\end{eqnarray}

Define $w=\frac{D^2f(0)(z^2)}{2!}$, then
\begin{eqnarray*}
D^2f(0)\Big(\frac{z+w}{2},\frac{z+w}{2}\Big)=\frac{1}{4}D^2f(0)(z^2)+\frac{1}{2}D^2f(0)(z,w)+\frac{1}{4}D^2f(0)(w^2),
\end{eqnarray*}
\begin{eqnarray*}
D^2f(0)\Big(\frac{z-w}{2},\frac{z-w}{2}\Big)=\frac{1}{4}D^2f(0)(z^2)-\frac{1}{2}D^2f(0)(z,w)+\frac{1}{4}D^2f(0)(w^2).
\end{eqnarray*}

Because  $D^2f_j(0)(z^2)=z_j\Big(\sum_{l=1}^{n}a_{jl}^{(j)}z_l\Big)$, where $|z_j|=\Arrowvert z\Arrowvert=\max_{1\le k\le n}\{|z_k|\}$, then
\begin{eqnarray}
D^2f_j(0)(z,w)&=&D^2f_j(0)\Big(\frac{z+w}{2},\frac{z+w}{2}\Big)-D^2f_j(0)\Big(\frac{z-w}{2},\frac{z-w}{2}\Big)\nonumber \\
&=&-\frac{1}{4}(z_j-w_j)\Big(\sum_{l=1}^{n}a_{jl}^{(j)}(z_l-w_l)\Big)\nonumber \\
&=&\frac{1}{2}\Big{[}z_j\Big(\sum_{l=1}^{n}a_{jl}^{j}w_l\Big)+w_j\Big(\sum_{l=1}^{n}a_{jl}^{(j)}z_l\Big)\Big{]}.
\end{eqnarray}

By (3.5) and (3.6), we have
\begin{eqnarray*}
\Big|D^2f_j(0)\Big(z_0,\frac{D^2f(0)(z_{0}^2)}{2!}\Big)\Big|&=&\frac{1}{2}\Big|\frac{z_j}{\Arrowvert z\Arrowvert}
\Big(\sum_{l=1}^{n}a_{jl}^{(j)}\frac{D^2f_l(0)(z_{0}^2)}{2!}\Big)+\frac{D^2f_j(0)(z_{0}^2)}{2!}\Big(\sum_{l=1}^{n}a_{jl}^{(j)}
\frac{z_l}{\Arrowvert z\Arrowvert}\Big)\Big|\\
&\le&\frac{1}{2}\Big[\Big|\sum_{l=1}^{n}a_{jl}^{(j)}\frac{D^2f_l(0)(z_{0}^2)}{2!}\Big|+\frac{|D^2f_j(0)(z_{0}^2)|}{2!}
\Big|\sum_{l=1}^{n}a_{jl}^{(j)}\frac{z_l}{\Arrowvert z\Arrowvert}\Big|\Big]\\
&\le&\frac{1}{2}\Big[\sum_{l=1}^{n}|a_{jl}^{(j)}|\frac{|D^2f_l(0)(z_{0}^2)|}{2!}+\frac{|D^2f_j(0)(z_{0}^2)|}{2!}
\sum_{l=1}^{n}|a_{jl}^{(j)}|\frac{|z_l|}{\Arrowvert z\Arrowvert}\Big]\\
&\le&\frac{1}{2}\Big[2(1-\alpha)(1-\alpha)+(1-\alpha)2(1-\alpha)\Big]\\
&=&2(1-\alpha)^2,j=1.2,\cdots,n.
\end{eqnarray*}
That is
\begin{eqnarray*}
\Big{\Arrowvert} D^2f(0)\Big(z_0,\frac{D^2f(0)(z_{0}^2)}{2!}\Big)\Big{\Arrowvert}\le2(1-\alpha)^2,\ z_0\in\partial D^n.
\end{eqnarray*}

In view of the maximum modulus theorem of holomorphic functions on the unit polydisk, we obtain
\begin{eqnarray*}
\Big{\Arrowvert} D^2f(0)\Big(z,\frac{D^2f(0)(z^2)}{2!}\Big)\Big{\Arrowvert}\le2(1-\alpha)^2\Arrowvert z\Arrowvert^3, z\in D^n.
\end{eqnarray*}
Using Corollary~2.2, we conclude that
\begin{eqnarray*}
\frac{\Arrowvert D^3f(0)(z^3)\Arrowvert}{3!}&\le&\frac{1-\alpha}{3}\Arrowvert z\Arrowvert^3+\frac{1}{3}\Big{\Arrowvert} D^2f(0)\Big(z,\frac{D^2f(0)(z^2)}{2!}\Big)\Big{\Arrowvert}\\
&\le&\frac{1-\alpha}{3}\Arrowvert z\Arrowvert^3+\frac{2(1-\alpha)^2}{3}\Arrowvert z\Arrowvert^3\\
&=&\frac{(1-\alpha)(3-2\alpha)}{3}\Arrowvert z\Arrowvert^3,\ z\in D^n.
\end{eqnarray*}

Finally, it is easy to prove the function
\begin{eqnarray*}
f(z)=\left\{
\begin{array}{lll}
\Big(\frac{1-(1-z_1)^{2\alpha-1}}{2\alpha-1},\frac{1-(1-z_2)^{2\alpha-1}}{2\alpha-1},\cdots,
\frac{1-(1-z_n)^{2\alpha-1}}{2\alpha-1}\Big)', & \alpha\in[0,1)\backslash\{\frac{1}{2}\},\\
(-\log(1-z_1),-\log(1-z_2),\cdots,-\log(1-z_n))',& \alpha=\frac{1}{2}
\end{array}
\right.
\end{eqnarray*}
meets the conditions of Theorem 3.2. Taking $z=(r,0,\cdots,0)'(0\le r<1)$, then
\begin{eqnarray*}
\frac{\Arrowvert D^3f(0)(z^3)\Arrowvert}{3!}=\frac{(1-\alpha)(3-2\alpha)}{3}r^3.
\end{eqnarray*}
This completes the proof of Theorem 3.2.

{\bf Theorem~3.3}\quad Suppose $0\le\alpha\le \frac{1}{2}$,  $f$ is an almost starlike of order $\alpha$ in $D^n$, and
\begin{eqnarray*}
&&2z_k\frac{D^2f_k(0)}{2!}\Big{(}z_0,\frac{D^2f(0)(z_0^2)}{2!}\Big{)}\|z\|-\frac{1}{1-\alpha}\Big{[}\frac{D^2f_k(0)(z_{0}^2)}{2!}\Big{]}^2\|z\|^2\\
&&=\frac{1-2\alpha}{1-\alpha}z_{k}^2\Big{(}\sum_{k=1}^na_{km}\frac{z_k}{\|z\|}(\sum_{l=1}^na_{ml}\frac{z_l}{\|z\|})\Big{)},\, \mbox{ for }\, z\in D^n\backslash\{0\}.
\end{eqnarray*}
where $a_{ml}=\frac{1}{2!}\frac{\partial^2f_m(0)}{\partial z_m\partial z_l}, m,l=1,2,\cdots,n$, then
\begin{eqnarray*}
\frac{\|D^3f(0)(z^3)\|}{3!}\le (1-\alpha)(3-4\alpha)\|z\|^3,\, \mbox{ for }\,  z\in D^n.
\end{eqnarray*}
The above estimate is sharp.

{\bf Proof.}\quad For fixed $z\in D^n\setminus\{0\}$, define $z_0=\frac{z}{\Arrowvert z\Arrowvert}$. Taking $T_z=(0,\cdots,0,\frac{|z_j|}{z_j},\cdots,0)$, where $|z_j|=\|z\|=\max_{1\le k\le n}\{|z_k|\}$. Applying Theorem~2.3, we have
\begin{eqnarray*}
\Big{|}2\frac{D^3f_j(0)(z_{0}^3)}{3!}\frac{\|z\|}{z_j}-2\frac{D^2f_j(0)}{2!}\Big{(}z_0,\frac{D^2f(0)(z_{0}^2)}{2!}\Big{)}\frac{\|z\|}{z_j}+\frac{1}{1-\alpha}\Big{[}\frac{D^2f_j(0)(z_{0}^2)}{2!}\Big{]}^2(\frac{\|z\|}{z_j})^2\Big{|}\le2(1-\alpha).
\end{eqnarray*}

According to the hypothesis of Theorem~3.3, we obtain
\begin{eqnarray*}
\frac{|D^3f_j(0)(z_{0}^3)|}{3!}\le 1-\alpha+\frac{1-2\alpha}{1-\alpha}\frac{M^2}{2},
\end{eqnarray*}
where $|z_j|=\|z\|=\max_{1\le k\le n}\{|z_k|\},\  M=\max_{1\le k\le n}\{\sum_{l=1}^n|a_{kl}|\}$. Applying the maximum modulus principle of holomorphic functions, and Lemmas~3.1, 3.3 and 3.4, $C_0=2(1-\alpha)$, then
\begin{eqnarray*}
\frac{\|D^3f(0)(z^3)\|}{3!}\le (1-\alpha+\frac{1-2\alpha}{1-\alpha}\frac{M^2}{2})\|z\|^3\le(1-\alpha)(3-4\alpha)\|z\|^3, z\in D^n.
\end{eqnarray*}

Finally, it is easy to prove that the function
\begin{eqnarray*}
f(z)=\left\{
\begin{array}{lll}
\Big{(}\frac{z_1}{(1-(1-2\alpha)z_1)^{\frac{2(1-\alpha)}{1-2\alpha}}},\cdots,\frac{z_n}{(1-(1-2\alpha)z_n)^{\frac{2(1-\alpha)}{1-2\alpha}}}\Big{)}',&0\leq\alpha<\frac{1}{2},\\
(z_1e^{z_1},\cdots,z_ne^{z_n})',& \alpha=\frac{1}{2}
\end{array}
\right.
\end{eqnarray*}
meets the hypothesis of Theorem~3.3, where $z\in D^n$. Taking $z=(r,0,\cdots,0)'(0\le r<1)$,
\begin{eqnarray*}
D^3f(0)(z^3)=(6(1-\alpha)(3-4\alpha)r^3,0,\cdots,0)'.
\end{eqnarray*}
Then
\begin{eqnarray*}
\frac{\|D^3f(0)(z^3)\|}{3!}=(1-\alpha)(3-4\alpha)r^3,
\end{eqnarray*}
This completes the proof of Theorem 3.3.

{\bf Corollary~3.3}\quad Suppose $0\le\alpha\le\frac{1}{2}$, $f$ is an almost starlike mapping of order $\alpha$ in $D^n$, and $\frac{D^2f_k(0)(z^2)}{2!}=z_k(\sum_{l=1}^n{a_{kl}z_l}), k=1,2,\cdots$, where $a_{ml}=\frac{1}{2!}\frac{\partial^2f_m(0)}{\partial z_m\partial z_l}, m,l=1,2,\cdots,n$, then
\begin{eqnarray*}
\frac{\|D^3f(0)(z^3)\|}{3!}\le (1-\alpha+\frac{1-2\alpha}{1-\alpha}\frac{M^2}{2})\|z\|^3\le(1-\alpha)(3-4\alpha)\|z\|^3, z\in D^n.
\end{eqnarray*}
where $M=\max_{1\le k\le n}\{\sum_{l=1}^n|a_{kl}|\}$. The above estimate is sharp.

{\bf Remark}\quad Setting $\alpha=0$ and $f\in S^*(D^n)$ in Corollary~3.3, we can obtain Theorem~1.3 in \cite{LL2012}.

{\bf Theorem~3.4}\quad  If $0\le\alpha<1$, $f$ is an almost starlike of order $\alpha$ in $D^n$, then
\begin{eqnarray*}
&&\frac{\|D^3f(0)(z^3)\|}{3!}\le(1-\alpha)(5-4\alpha)\|z\|^3,z\in D^n,n=2,\\
&&\frac{\|D^3f(0)(z^3)\|}{3!}\le(1-\alpha)(3\sqrt{3}+1-3\sqrt{3}\alpha)\|z\|^3,z\in D^n,n\ge3,n\in\mathbb{N}_+.
\end{eqnarray*}

{\bf Proof}\quad For fixed $z\in D^n\setminus\{0\}$, define $z_0=\frac{z}{\Arrowvert z\Arrowvert}$. By Lemmas~2.3 and ~3.1, we have
\begin{eqnarray}
\Arrowvert D^2f(0)(z^{2}_0)\Arrowvert\le 4(1-\alpha),\|z_0\|=1.
\end{eqnarray}

Apparently, $\hat{L}(z^2)=D^2f(0)(z^{2})$ is a bounded linear operator, we have
\begin{eqnarray}
\|\hat{L}\|=\sup{\{\|\hat{L}(z^2)\|:\|z\|\le1\}}=\sup{\{\|\hat{L}(z^2)\|:\|z\|=1\}}.
\end{eqnarray}

According to (3.3),\ (3.4), we obtain
\begin{eqnarray*}
\|\hat{L}\|=4(1-\alpha).
\end{eqnarray*}

Let $L(x,y)=D^2f(0)(x,y),x,y\in X$, by Lemma~3.6 we have
\begin{eqnarray*}
\|L\|&=&2(1-\alpha),\, \mbox{ for }\, n=2,\\
\|L\|&\le&\frac{3}{2}\sqrt{3}(1-\alpha),\, \mbox{ for }\, n\ge3,n\in\mathbb{N}_+.
\end{eqnarray*}

Let $w=\frac{D^2f(0)(z^2)}{2!}$ and $\| w_0\|=\frac{\| D^2f(0)(z_{0}^2)\|}{2!}\le 2(1-\alpha)$, again by Lemma~3.7,
\begin{eqnarray*}
\|D^2f(0)(z,w)\|&\le &\|L\|\|z\|\|w\|\le4(1-\alpha)\|z\|\|w\|\le8(1-\alpha)^2\|z\|^2, z\in D^n,n=2,\\
\|D^2f(0)(z,w)\|&\le &\|L\|\|z\|\|w\|\le3\sqrt{3}(1-\alpha)\|z\|\|w\|\\
&\le &6\sqrt{3}(1-\alpha)^2\|z\|^2, z\in D^n,n\ge3,n\in\mathbb{N}_+.
\end{eqnarray*}

According to Lemmas~2.2 and 2.4, we conclude that
\begin{eqnarray}
\Big{|}2\frac{D^3f(0)(z_0^3)}{3!}-2\frac{D^2f(0)}{2!}(z_0,w_0)\Big{|}\le2(1-\alpha).
\end{eqnarray}

By (3.7), (3.8) and (3.9), we have
\begin{eqnarray*}
&&\frac{\Arrowvert D^3f(0)(z_0^3)\Arrowvert}{3!}\le {1-\alpha+\frac{1}{2}\|Df(0)(z_0,w_0)\|}\le(1-\alpha)(5-4\alpha),n=2,\\
&&\frac{\Arrowvert D^3f(0)(z_0^3)\Arrowvert}{3!}\le(1-\alpha)(3\sqrt{3}+1-3\sqrt{3}\alpha),n\ge3,n\in\mathbb{N}_+,
\end{eqnarray*}
that is
\begin{eqnarray*}
&&\frac{\|D^3f(0)(z^3)\|}{3!}\le(1-\alpha)(5-4\alpha)\|z\|^3,z\in D^n,n=2,\\
&&\frac{\|D^3f(0)(z^3)\|}{3!}\le(1-\alpha)(3\sqrt{3}+1-3\sqrt{3}\alpha)\|z\|^3,z\in D^n,n\ge3,n\in\mathbb{N}_+.
\end{eqnarray*}
This completes the proof.

{\bf Theorem~3.5}\quad Suppose $0\le\alpha\le\frac{37-\sqrt{505}}{36}$, $f$ is an almost starlike mapping of order $\alpha$ in $D^n$, and $\frac{D^2f_k(0)(z^2)}{2!}=a_kz_{k}^2, \frac{D^3f_k(0)(z^3)}{3!}=b_kz_{k}^3, k=1,2,\cdots,n, z\in D^n$, where
$a_{k}=\frac{1}{2!}\frac{\partial^2f_k(0)}{\partial^2 z_k}, b_{kl}=\frac{1}{3!}\frac{\partial^3f_k(0)}{\partial^3 z_k}, k,l=1,2,\cdots,n$, then
\begin{eqnarray*}
\frac{\|D^4f(0)(z^4)\|}{4!}\le \frac{(1-\alpha)(3-4\alpha)(4-6\alpha)}{3}\|z\|^4, z\in D^n.
\end{eqnarray*}
The above estimate is sharp.

{\bf Proof}\quad For fixed $z\in D^n\backslash\{0\}$, let $z_0=\frac{z}{\|z\|}$. Define $p_j(\xi)=\frac{g_j(\xi z_0)\|z\|}{\xi z_j}, \xi\in D$, where
$$g(z)=(g_1(z),g_2(z),\cdots,g_n(z))'=(Df(z))^{-1}f(z)$$ is a column vector in $\mathbb{C}^n$, $|z_j|=\|z\|=\max_{1\le k\le n}\{|z_k|\}$. Since $f(z)$ is an almost starlike mapping of order $\alpha$ in $D^n$, we have $p_j(\xi)\in H(D), {\rm Re}p_j(\xi)\ge \alpha, p_j(0)=1$, and
\begin{eqnarray*}
p_j(\xi)=1+\frac{D^2g_j(0)(z_{0}^2)}{2!}\frac{\|z\|}{z_j}\xi+\cdots+\frac{D^{m}g_j(0)(z_{0}^{m})}{m!}\frac{\|z\|}{z_j}\xi^{m-1}+\cdots,\xi\in D.
\end{eqnarray*}

According to Lemma~2.2, we have
\begin{eqnarray}
&&\Big{|}\frac{D^{4}g_j(0)(z_{0}^{4})}{4!}\frac{\|z\|}{z_j}-\frac{\frac{D^{2}g_j(0)(z_{0}^{2})}{2!}\frac{\|z\|}{z_j}\frac{D^{3}g_j(0)(z_{0}^{3})}{3!}
\frac{\|z\|}{z_j}}{1-\alpha}
+\frac{(\frac{D^{2}g_j(0)(z_{0}^{2})}{2!}\frac{\|z\|}{z_j})^3}{4(1-\alpha)^2}\Big{|}\nonumber\\
&\le& 2(1-\alpha)-\frac{|\frac{D^{2}g_j(0)(z_{0}^{2})}{2!}\frac{\|z\|}{z_j}|^2}{2(1-\alpha)}.
\end{eqnarray}

From Lemma~2.4 and the conditions of Theorem~3.4, we have
\begin{eqnarray*}
&&\frac{D^{2}g_j(0)(z_{0}^{2})}{2!}\frac{\|z\|}{z_j}=-a_j\frac{z_j}{\|z\|},\\
&&\frac{D^{3}g_j(0)(z_{0}^{3})}{3!}\frac{\|z\|}{z_j}=2a_{j}^2(\frac{z_j}{\|z\|})^2-2b_j(\frac{z_j}{\|z\|})^2,\\
&&\frac{D^{4}g_j(0)(z_{0}^{4})}{4!}\frac{\|z\|}{z_j}=-3\frac{D^{4}f_j(0)(z_{0}^{4})}{4!}\frac{\|z\|}{z_j}-4a_{j}^3(\frac{z_j}{\|z\|})^3
+7a_{j}b_{j}(\frac{z_j}{\|z\|})^3.
\end{eqnarray*}

Consequently, we have
\begin{eqnarray}
-3\frac{D^{4}f_j(0)(z_{0}^{4})}{4!}\frac{\|z\|}{z_j}=\frac{D^{4}g_j(0)(z_{0}^{4})}{4!}\frac{\|z\|}{z_j}+4a_{j}^3(\frac{z_j}{\|z\|})^3
-7a_{j}b_{j}(\frac{z_j}{\|z\|})^3,
\end{eqnarray}
and
\begin{eqnarray}
\Big{|}\frac{D^{4}g_j(0)(z_{0}^{4})}{4!}\frac{\|z\|}{z_j}+\frac{7-8\alpha}{4(1-\alpha)^2}a_{j}^3(\frac{z_j}{\|z\|})^3-
\frac{2}{1-\alpha}a_{j}b_{j}(\frac{z_j}{\|z\|})^3\Big{|}\le 2(1-\alpha)-\frac{|a_j|^2}{2(1-\alpha)}.
\end{eqnarray}

According to Theorem~2.3, we have
\begin{eqnarray}
\Big{|}\frac{3-4\alpha}{2(1-\alpha)}a_{j}^2(\frac{z_j}{\|z\|})^2-2b_j(\frac{z_j}{\|z\|})^2\Big{|}\le2(1-\alpha)-\frac{|a_j|^2}{2(1-\alpha)}.
\end{eqnarray}

Connecting (3.11), (3.12) and (3.13), we have
\begin{eqnarray*}
\Big{|}-3\frac{D^{4}f_j(0)(z_{0}^{4})}{4!}\frac{\|z\|}{z_j}\Big{|}&=&\Big{|}\frac{D^{4}g_j(0)(z_{0}^{4})}{4!}\frac{\|z\|}{z_j}+
4a_{j}^3(\frac{z_j}{\|z\|})^3-7a_{j}b_{j}(\frac{z_j}{\|z\|})^3\Big{|}\nonumber\\
&=&\Big{|}\Big{(}\frac{D^{4}g_j(0)(z_{0}^{4})}{4!}\frac{\|z\|}{z_j}+\frac{7-8\alpha}{4(1-\alpha)^2}a_{j}^3(\frac{z_j}
{\|z\|})^3-\frac{2}{1-\alpha}a_{j}b_{j}(\frac{z_j}{\|z\|})^3\Big{|}\Big{)}\nonumber\\
&&+\frac{5-7\alpha}{2(1-\alpha)}a_j\frac{z_j}{\|z\|}\Big{[}\frac{3-4\alpha}{2(1-\alpha)}a_{j}^2(\frac{z_j}{\|z\|})^2
-2b_j(\frac{z_j}{\|z\|})^2\Big{]}\\
&&-\frac{12\alpha^2-17\alpha+6}{4(1-\alpha)^2}a_{j}^3(\frac{z_j}{\|z\|})^3\Big{|}\nonumber\\
&=&2(1-\alpha)-\frac{|a_j|^2}{2(1-\alpha)}+\frac{5-7\alpha}{2(1-\alpha)}|a_j|\Big{(}2(1-\alpha)
-\frac{|a_j|^2}{2(1-\alpha)}\Big{)}\nonumber\\
&&+\frac{12\alpha^2-17\alpha+6}{4(1-\alpha)^2}|a_{j}|^3\nonumber\\
&=&\frac{12\alpha^2-10\alpha+1}{4(1-\alpha)^2}|a_{j}|^3-\frac{|a_j|^2}{2(1-\alpha)}+(5-7\alpha)|a_j|+2(1-\alpha).
\end{eqnarray*}
By Lemma~3.8, we obtain
\begin{eqnarray*}
\Big{|}\frac{D^{4}f_j(0)(z_{0}^{4})}{4!}\Big{|}&\leq &\frac{1}{3}\Big{[}\frac{12\alpha^2-10\alpha+1}{4(1-\alpha)^2}M^3
-\frac{1}{2(1-\alpha)}M^2+(5-7\alpha)M|+2(1-\alpha)\Big{]}\\
&\leq &\frac{(1-\alpha)(3-4\alpha)(4-6\alpha)}{3}.
\end{eqnarray*}

Notice that $M=\max_{\{1\le k\le n\}}\{|a_k|\}\le2(1-\alpha)$, by Lemma~3.4,
\begin{eqnarray*}
\frac{\|D^{4}f(0)(z^{4})\|}{4!}\le\frac{(1-\alpha)(3-4\alpha)(4-6\alpha)}{3}\|z\|^4, z\in D^n.
\end{eqnarray*}

Example in the proof of Theorem 3.3 verifies the accuracy of theorem 3.5. This completes the proof.

\section*{Acknowledgment}
The authors thank Professor Tai-Shun Liu for providing a lot of valuable suggestions to improve the quality of this paper.

\end{document}